\documentclass[oneside,english]{amsart}
\usepackage[T1]{fontenc}
\usepackage[latin9]{inputenc}
\usepackage{amsthm}
\usepackage[pdftex]{graphicx}

\makeatletter
\numberwithin{equation}{section}
\numberwithin{figure}{section}

\makeatother

\usepackage{babel}
\begin{document}

\title[\tiny{Boundary Value Problems for a Family of Domains in the Sierpinski
Gasket}]{Boundary Value Problems for a Family of Domains in the Sierpinski
Gasket}
%    author one information
\author{Zijian Guo}
\address{Department of Mathematics, The Chinese University of Hong Kong, Shatin, Hong Kong.}
\curraddr{Department of Statistics, The Wharton School, University
of Pennsylvania, U.S.A.} \email{zijguo@wharton.upenn.edu}
\thanks{}

%    author two information
\author{Hua Qiu}
\address{Department of Mathematics, Nanjing University, Nanjing, 210093, P. R. China.}
\curraddr{} \email{huatony@gmail.com}
\thanks{The research of the second author was supported by the Nature Science Foundation of Jiangsu Province of China, Grant
BK20131265, and the Project Funded by the Priority Academic Program
Development of Jiangsu Higher Education Institutions.}

\author{Robert S. Strichartz}
\address{Department of Mathematics, Malott Hall, Cornell University, Ithaca, NY 14853, U.S.A.}
\curraddr{} \email{str@math.cornell.edu}
\thanks{The research of the third author was supported in part by the National Science foundation, grant DMS-1162045.}

%    For articles to be published after 1 January 2010, you may use
%    the following version:
%\subjclass[2010]{Primary }


\date{}

\dedicatory{}
\begin{abstract}
For a family of domains in the Sierpinski gasket, we study harmonic
functions of finite energy, characterizing them in terms of their
boundary values, and study their normal derivatives on the boundary.
We characterize those domains for which there is an extension operator
for functions of finite energy. We give an explicit construction of
the Green's function for these domains.
\end{abstract}
\maketitle

\section{Introduction }

Consider the domain $\Omega_{x}$ in the Sierpinski Gasket ($\mathcal{SG}$) consisting
of all points above the horizontal line $L_{x}$ at the distance $x$
from the top vertex $q_{0}$, for $0<x\leq1$.\\

\begin{figure}[h]
\begin{center}
\includegraphics[width=7.8cm,totalheight=4.6cm]{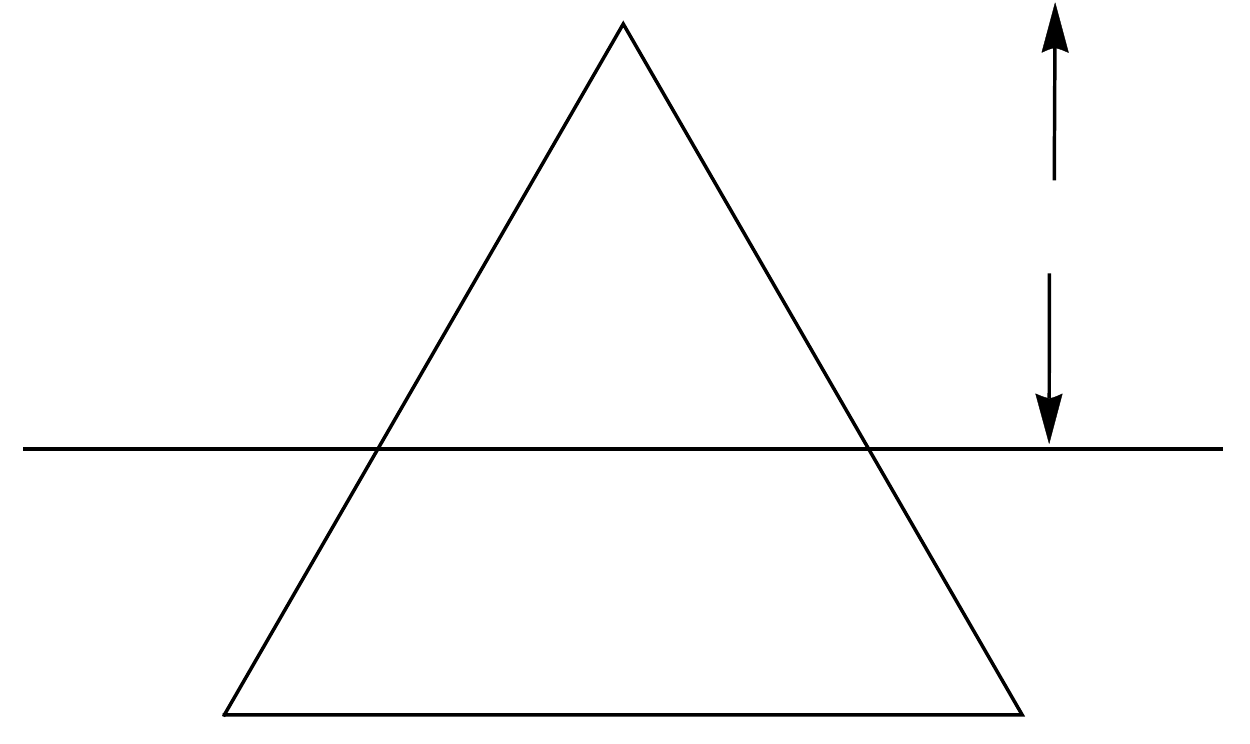}
\setlength{\unitlength}{1cm}
\begin{picture}(0,0) \thicklines
\put(-1.55,3.1){$x$}
\put(-4.3,4.7){$q_{0}$}
\put(-4.3,2.7){$\Omega_x$}
\put(-.3,1.8){$L_x$}
\end{picture}
\begin{center}
\textbf{Figure 1.1.}
\end{center}
\end{center}
\end{figure}

Let $S(x)=\mathcal{SG}\cap L_{x}$. For $x$ not a dyadic rational, this is
a Cantor set. The boundary of $\Omega_{x}$ consists of $S(x)$ together
with $q_{0}$. By general principles, harmonic functions on $\Omega_{x}$
are determined by their boundary values, where harmonic functions
are defined to be solutions of $\triangle h=0$ on the interior of
$\Omega_{x}$, where $\triangle$ is the Kigami Laplacian on $\mathcal{SG}$.
The study of such harmonic functions was initiated in {[}S1{]}, and
continued in {[}OS{]} for the special case $x=1$. In this paper,
we extend the results in {[}OS{]} to the general case. In Section
2, we give an explicit description of the analog of the Poisson kernel
to recover the harmonic function from its boundary values, in terms
of the Haar series expansion of the boundary values on $S(x)$, and
we characterize the boundary values that correspond to harmonic functions
of finite energy. In Section 3, we define normal derivatives on the
boundary and give a description of the Dirichlet-to-Neumann map
as a multiplier transform on the Haar series expansion.

In Section 4, we study the extension problem for functions of finite
energy on $\Omega_{x}$ to functions of finite energy on
$\mathcal{SG}$. We are able to characterize the values of $x$ for
which such extensions are possible. In particular, the value $x=1$
studied in {[}OS{]} does not admit such extensions. This may be
regarded as the first of a family of Sobolev extension problems,
based on Sobolev spaces on $\mathcal{SG}$ discussed in {[}S2{]}. We
leave these as open problems for future research. Related problems
are studied in {[}LS{]} and {[}LRSU{]}.

In Section 5, we give a construction of a Green's function on $\Omega_{x}$
to solve the Dirichlet problem $-\triangle u=F$ on $\Omega_{x}$,
$u|_{\partial\Omega_{x}}=0$ via an integral transform of $F$. The
construction of the Green's function is analogous to Kigami's construction
on $\mathcal{SG}$.

The reader is referred to the books {[}Ki{]} and {[}S3{]} for a description
of the theory of the Laplacian on $\mathcal{SG}$, and related fractals. It
would be interesting to extend the results of this paper to other
domains in $\mathcal{SG}$, and to domains in other fractals. In this regard, we
offer the following cautionary tale. Consider
the fractal $\mathcal{SG}_{3}$, defined similarly to $\mathcal{SG}$ but by subdivisions
of the sides of triangles into three rather than two pieces(see Figure
1.2).

\begin{figure}[h]
\begin{center}
\includegraphics[width=7.8cm,totalheight=7.8cm]{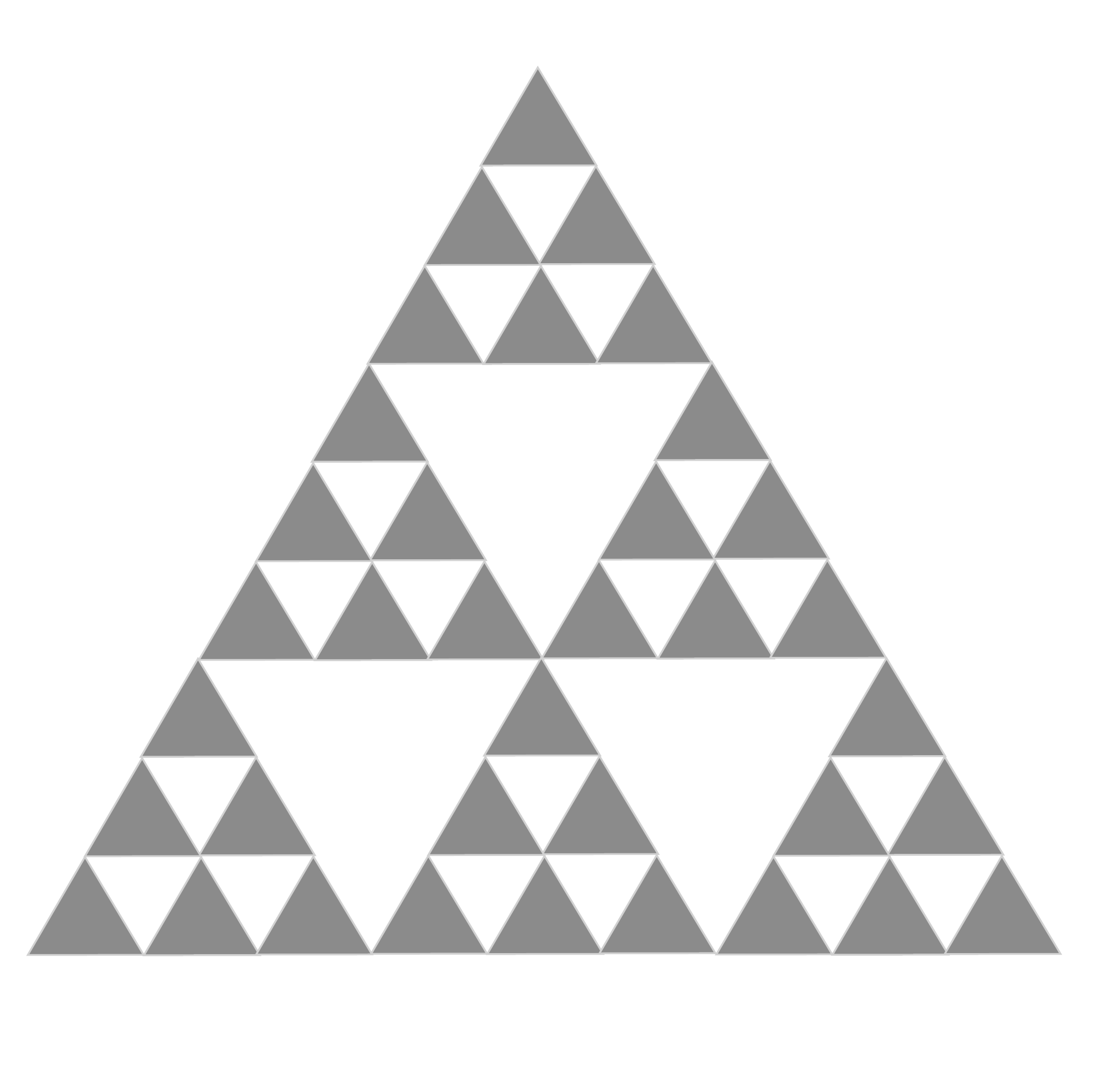}
\setlength{\unitlength}{1cm}
\begin{center}
\begin{picture}(0,0) \thicklines
\put(0.7,0.9){$\mathcal{SG}_3$}
\end{picture}
\textbf{Figure 1.2.}
\end{center}
\end{center}
\end{figure}

We may consider domains $\Omega_{x}$ defined as before, with the
boundary $S(x)$ modeled as a Cantor set with divisions into three
pieces. There is a natural analog of Haar functions on $S(x)$, with
two generators as shown in Figure 1.3.

\begin{figure}[h]
\begin{center}
\includegraphics[width=7.8cm,totalheight=2.6cm]{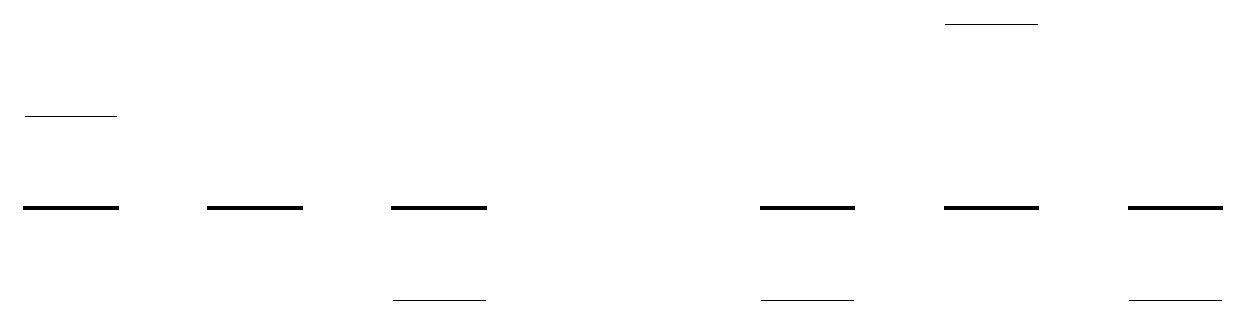}
\setlength{\unitlength}{1cm}
\begin{picture}(0,0) \thicklines
\put(-8.3,1.6){$1$}
\put(-8.5,0.1){$-1$}
\put(-4,0.1){$-1$}
\put(-3.8,2.3){$2$}
\end{picture}\\
\begin{center}
\textbf{\\Figure 1.3.}\small{ Haar generators.}
\end{center}
\end{center}
\end{figure}

Because the second generator is symmetric rather than skew-symmetric,
we cannot glue to zero at the top, so the analog of Lemma 2.3 does
not hold. It is not clear how to overcome this difficulty.

\section{Harmonic Functions on $\Omega_{x}$}

For $0<x\leq1$, there is a unique representation

\begin{equation}
x=\sum_{k=1}^{\infty}2^{-n_{k}}\label{eq:2.1}
\end{equation}
for a sequence

\begin{equation}
0<n_{1}<n_{2}<\cdots\label{eq:2.2}
\end{equation}
of increasing positive integers. We will approximate $\Omega_{x}$ by
the increasing sequence of domains $\Omega_{x}^{(m)}$ where each
$\Omega_{x}^{(m)}$ is the closure of $\Omega_{x_{[m]}}$ where
\begin{equation}
x_{[m]}=\sum_{k=1}^{m}2^{-n_{k}}\label{eq:2.3}
\end{equation}
is the partial sum of $\eqref{eq:2.1}$. (Note that $\eqref{eq:2.3}$
is not the representation of $x_{[m]}$ of the form $\eqref{eq:2.1}$
since it is a finite binary representation.) The domain
$\Omega_x^{(m)}$ is a finite union of cells, specifically $1$
$n_1$-cell, $2$ $n_2$-cells, $4$ $n_3$-cells, $\cdots$, $2^{m-1}$
$n_m$-cells.  Figure 2.1 illustrates $\Omega_x^{(m)}$ for $m=1,2,3$
for two choices of $x$. The boundary of $\Omega_x^{(m)}$ consists of
the top vertex $q_0$ together with the $2^m$ bottom vertices of the
$n_m$-cells.

\begin{figure}[h]
\begin{center}
\includegraphics[width=12cm,totalheight=3.9cm]{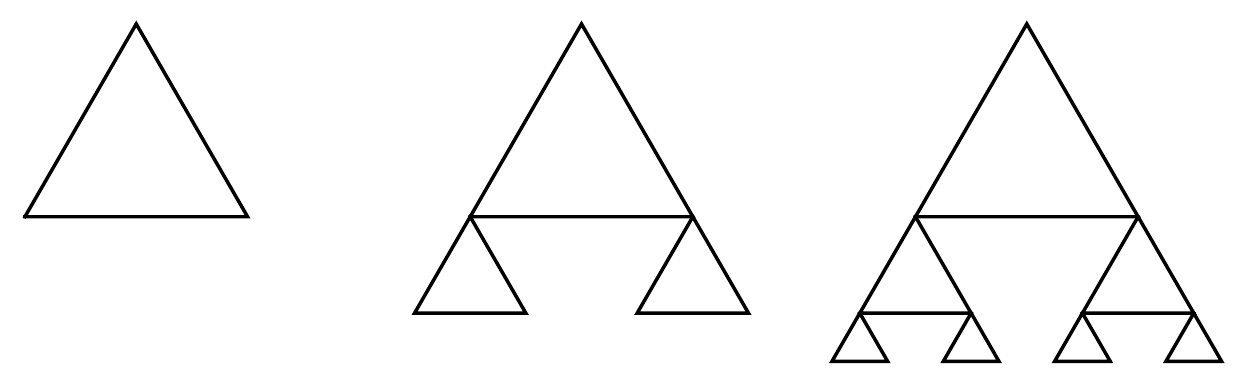}\vspace{1.2cm}
\includegraphics[width=12cm,totalheight=3.2cm]{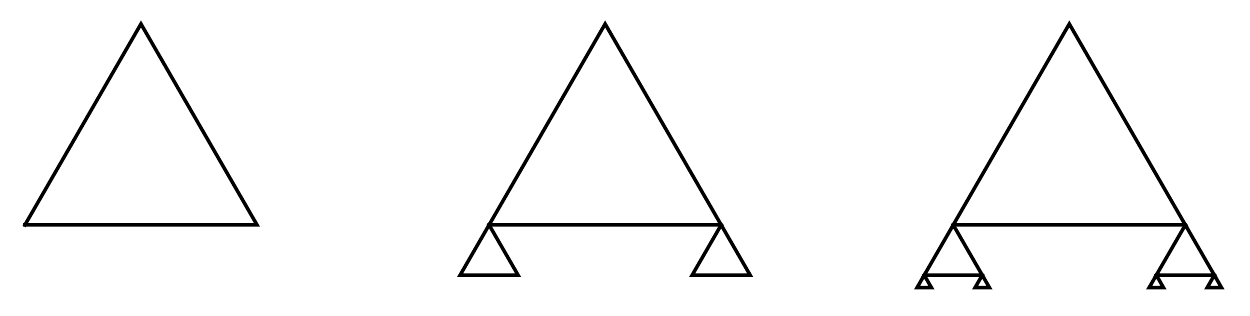}
\setlength{\unitlength}{1cm}
\begin{picture}(0,0) \thicklines
\put(-2.8,4.3){$m=3$} \put(-7.0,4.3){$m=2$} \put(-11.3,4.3){$m=1$}
\put(-8.0,3.8){$n_1=1, n_2=2, n_3=3$} \put(-2.5,-0.1){$m=3$}
\put(-6.9,-0.1){$m=2$} \put(-11.3,-0.1){$m=1$}
\put(-8.0,-0.6){$n_1=1, n_2=3, n_3=5$}
\end{picture}
\begin{center}\vspace{1cm}
\textbf{Figure 2.1.}\small{ Some examples of $\Omega_x^{(m)}$ for
$m=1,2,3$.}
\end{center}
\end{center}
\end{figure}

Following {[}S1{]} we define
\begin{equation}
Rx=\sum_{k=2}^\infty 2^{-n_k}=x-2^{-n_1}\label{eq:2.4}
\end{equation}
and the function $m_0(x)$ by the identity
\begin{equation}
m_0(x)=\frac{1}{1+2(\frac{5}{3})^{n_2-n_1}(1-m_0(Rx))}\label{eq:2.5}
\end{equation}
which is easily solved to obtain a variant of a continued fraction representation
\begin{equation}
m_0(x)=\lim_{k\rightarrow\infty}m_0^{(k)}(x)\label{eq:2.6}
\end{equation}
for
\begin{equation}
m_0^{(k)}=\dfrac{1}{1+2\left(\frac{5}{3}\right)^{n_2-n_1}(1-\dfrac{1}{1+2\left(\frac{5}{3}\right)^{n_3-n_2}(1-\dfrac{1}{\ddots\dfrac{1}{1+2\left(\frac{5}{3}\right)^{n_k-n_{k-1}}}}}}.\label{eq:2.7}
\end{equation}
See Figure 2.2 for the graph of $m_0(x)$ on $(0,1]$.

\begin{figure}[h]
\begin{center}
\includegraphics[width=12cm,totalheight=7.9cm]{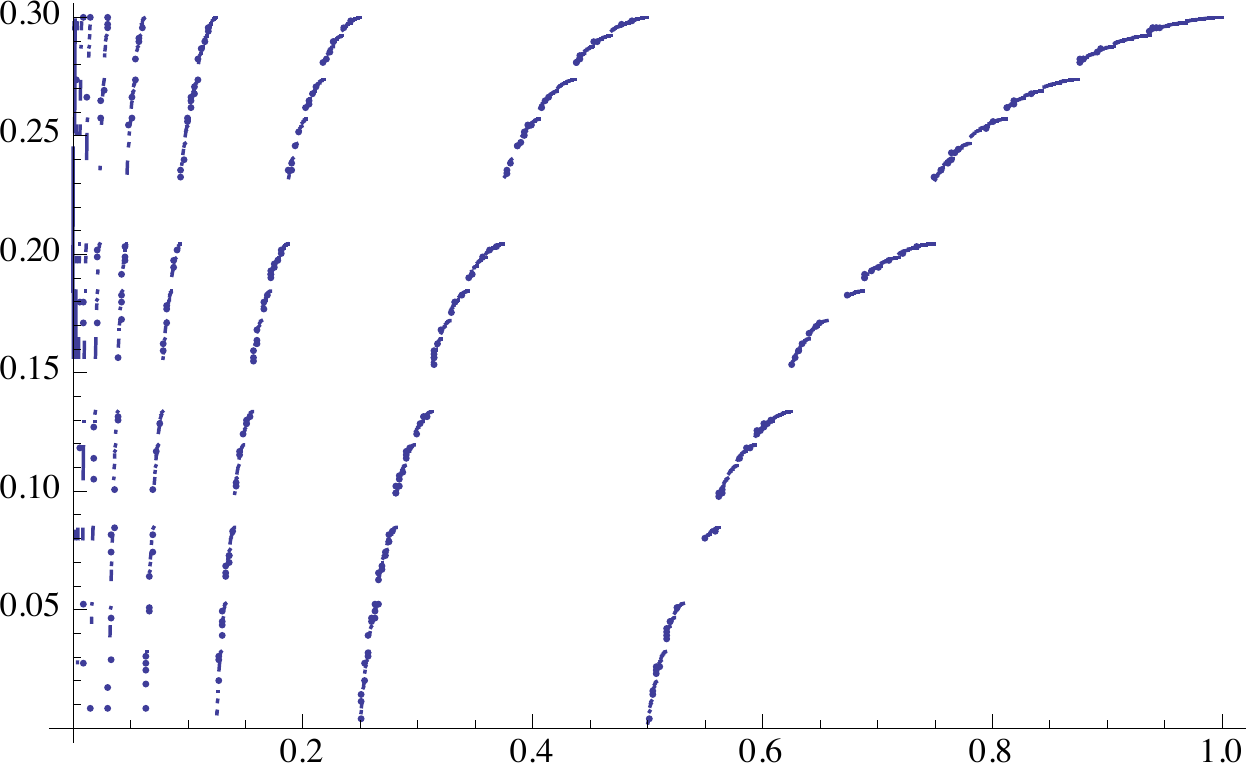}
\setlength{\unitlength}{1cm}
\begin{center}
\textbf{Figure 2.2.}\small{ The graph of $m_0(x)$.}
\end{center}
\end{center}
\end{figure}

 We also define
\begin{equation}
m_1(x)=\frac{1-m_0(x)^2}{2m_0(x)+1}, \quad m_2(x)=\frac{m_0(x)-m_0(x)^2}{2m_0(x)+1}.\label{eq:2.8}
\end{equation}

Note that
\begin{equation}
m_0(x)+m_1(x)+m_2(x)=1.\label{eq:2.9}
\end{equation}

These functions enable us to describe harmonic functions in $\Omega_x$. The boundary of $\Omega_x$ consists of the top vertex $q_0$ and $S(x)=L_x\cap \mathcal{SG}$. If $x$ is not a dyadic rational then $S(x)$ is a Cantor set. We will assume this holds. Then a harmonic function is determined by the value $h(q_0)$ and the expansion of $h|_{\mathcal{SG}}$ in a Haar basis.

\textbf{Definition 2.1. } The harmonic function $h_0$ satisfies
\begin{equation}
h_0(q_0)=1,\quad h_0|_{S(x)}=0.\label{eq:2.10}
\end{equation}
The harmonic function $h_1$ satisfies
\begin{equation}
h_1(q_0)=0,\quad h_1|_{S(x)\cap F_0^{n_1-1}F_1(\mathcal{SG})}=1,\quad h_1|_{S(x)\cap F_0^{n_1-1}F_2(\mathcal{SG})}=-1.\label{eq:2.11}
\end{equation}
We write $h_0^x$ and $h_1^x$ when we need to explicitly show the dependence on $x$.

Note that $1-h_0$ satisfies
\begin{equation}
(1-h_0)(q_0)=0,\quad (1-h_0)|_{S(x)}=1,
\end{equation}
so that $1-h_0$ and $h_1$ vanish at $q_0$ and give the first two
Haar functions when restricted to $S(x)$. Also it is shown in
{[}S1{]} that
\begin{equation}
h_0(F_0^{n_1-1}F_1q_0)=h_0(F_0^{n_1-1}F_2q_0)=m_0(x)\label{eq:2.13}
\end{equation}
and
\begin{equation}
h_1(F_0^{n_1-1}F_1q_0)=-h_1(F_0^{n_1-1}F_2q_0)=m_1(x)-m_2(x).\label{eq:2.14}
\end{equation}

\textbf{Lemma 2.2.} \textit{Let $y=2^{n_1}Rx$. Then
\begin{equation}
h_0^x\circ(F_0^{n_1-1}F_1)=h_0^x\circ(F_0^{n_1-1}F_2)=m_0(x)h_0^y\label{eq:2.15}
\end{equation}
and
\begin{equation}
h_1^x\circ(F_0^{n_1-1}F_1)=-h_1^x(F_0^{n_1-1}F_2)=1+(m_1(x)-m_2(x)-1)h_0^y.\label{eq:2.16}
\end{equation}}

\textit{Proof.} The function $m_0(x)h_0^y$ is a harmonic function on $\Omega_y$ with boundary values $m_0(x)$ at $q_0$ and zero on $S(y)$. Note that
$F_0^{n_1-1}F_1(S(y))=S(x)$, so $h_0^x\circ(F_0^{n_1-1}F_1)$ is also a harmonic function on $\Omega_y$ vanishing on $S(y)$, and it assume the value $m_0(x)$ at $q_0$ by $\eqref{eq:2.13}$. Thus $\eqref{eq:2.15}$ holds. A similar argument shows that $\eqref{eq:2.14}$ implies $\eqref{eq:2.16}.\Box$

Next we consider the general Haar basis functions on $L^2(S(x))$. Let $\omega=(\omega_1,\cdots,\omega_m)$ be a word of length $|\omega|=m$, with each $\omega_j=1$ or $2$. Then
\begin{equation}
S_\omega(x)=S(x)\cap F_0^{n_1-1}F_{\omega_1}F_0^{n_2-n_1-1}F_{\omega_2}\cdots F_0^{n_m-n_{m-1}-1}F_{\omega_m}(\mathcal{SG})\label{eq:2.17}
\end{equation}
describe the dyadic pieces of $S(x)$. In particular,
\begin{equation}
S(x)=\bigcup_{|\omega|=m}S_{\omega}(x).\label{eq:2.18}
\end{equation}
The Cantor measure $\mu$ on $S(x)$ assigns measure $2^{-m}$ to each piece $S_{\omega}(x)$. The
Haar function $\psi_{\omega}$ is supported on $S_{\omega}(x)$ and satisfies
\begin{equation}
\psi_\omega|_{S_{\omega1}(x)}=2^{m/2} \mbox{ and } \psi_\omega|_{S_{\omega2}(x)}=-2^{m/2}.\label{eq:2.19}
\end{equation}
Then $1\cup\{\psi_\omega\}$ is an orthonormal basis for $L^2(S(x),d\mu)$. We define $h_{\omega}^x$ to be the harmonic function on $\Omega_x$ with boundary values $h_\omega^x(q_0)=0$ and $h_{\omega}^x|_{S(x)}=\psi_\omega$.

\textbf{Lemma 2.3.} \textit{Let $y_m=2^{n_m}R^mx$. Then $h_\omega^x$ is supported in $$\Omega_x\cap F_0^{n_1-1}F_{\omega_1}F_{0}^{n_2-n_1-1}F_{\omega_2}\cdots F_0^{n_m-n_{m-1}-1}F_{\omega_m}(\mathcal{SG})$$and
\begin{equation}
h_\omega^x\circ(F_0^{n_1-1}F_{\omega_1}F_0^{n_2-n_1-1}F_{\omega_2}\cdots F_0^{n_m-n_{m-1}-1}F_{\omega_m})=2^{m/2}h_1^{y_m}.\label{eq:2.20}
\end{equation}}

\textit{Proof.} The key observation is that, because of skew-symmetry, the function $h_1$ not only vanishes at $q_0$ but also has normal derivative vanishing at $q_0$. Thus we may glue the function defined by $\eqref{eq:2.20}$ to zero outside this cell and still have a harmonic function. This function clearly has the required boundary values for $h_\omega^x. \Box$

\textbf{Theorem 2.4.} \textit{The energies are given by
\begin{equation}
\mathcal{E}(h_0^x)=(1-m_0(x))^2\sum_{j=1}^{\infty}2^{2-j}\left(\frac{5}{3}\right)^{2n_1-n_j},\label{eq:2.21}
\end{equation}
\begin{equation}
\mathcal{E}(h_1^x)=6\left(\frac{1-m_0(x)}{2m_0(x)+1}\right)^2\left(\frac{5}{3}\right)^{n_1}+2\left(\frac{3m_0(x)}{2m_0(x)+1}\right)^2\left(\frac{5}{3}\right)^{n_1}\mathcal{E}(h_0^y),\label{eq:2.22}
\end{equation}
and
\begin{equation}
\mathcal{E}(h_\omega^x)=2^m\left(\frac{5}{3}\right)^{n_m}\mathcal{E}(h_1^{y_m})\label{eq:2.23}
\end{equation}
where $m=|\omega|. $ Moreover, there exist positive constants $C_1$ and $C_2$, independent of $x$, such that
\begin{equation}
C_12^m\left(\frac{5}{3}\right)^{n_{m+1}}\leq \mathcal{E}(h_\omega^x)\leq C_2 2^m\left(\frac{5}{3}\right)^{n_{m+1}}.\label{eq:2.24}
\end{equation}}

\textit{Proof.} We compute the energy of $h_0^x$ on the top cell $F_0^{n_1}(\mathcal{SG})$ using $\eqref{eq:2.13}$ to be $\left(\frac{5}{3}\right)^{n_1}2(m_0(x)-1)^2$, since there are two edges where the difference of $h_0^x$ is $m_0(x)-1$. On the remaining cells $F_0^{n_1-1}F_1(\mathcal{SG})$ and $F_0^{n_1-1}F_2(\mathcal{SG})$ the function $h_0^x$ is equal to $m_0(x)h_0^y\circ(F_0^{n_1-1}F_1)^{-1}$ and $m_0(x)h_0^y\circ(F_0^{n_1-1}F_2)^{-1}$ by $\eqref{eq:2.15}$.  These each have energy $m_0(x)^2\left(\frac{5}{3}\right)^{n_1}\mathcal{E}(h_0^y)$, so
\begin{equation}
\mathcal{E}(h_0^x)=2\left(\frac{5}{3}\right)^{n_1}((m_0(x)-1)^2+m_0(x)^2\mathcal{E}(h_0^y)).\label{eq:2.25}
\end{equation}
Before iterating this identity we observe that
\begin{equation}
(1-m_0(y))m_0(x)=\frac{1}{2}\left(\frac{5}{3}\right)^{n_1-n_2}(1-m_0(x))\label{eq:2.26}
\end{equation}
as an immediate consequence of $\eqref{eq:2.5}$. Thus
$$\mathcal{E}(h_0^x)=(1-m_0(x))^2\left(2\left(\frac{5}{3}\right)^{n_1}+\left(\frac{5}{3}\right)^{2n_1-n_2}\right)+4\left(\frac{5}{3}\right)^{n_2}m_0(x)^2m_0(y)^2\mathcal{E}(h_0^{y_2})$$
and by iterating we obtain  $\eqref{eq:2.21}$.

Similarly, we use $\eqref{eq:2.14}$ to compute the energy of $h_1^x$ on the top cell $F_0^{n_1}(\mathcal{SG})$ to be $\left(\frac{5}{3}\right)^{n_1}6(m_2(x)-m_1(x))^2=\left(\frac{5}{3}\right)^{n_1}6\left(\frac{1-m_0(x)}{2m_0(x)+1}\right)^2$ by $\eqref{eq:2.8}$. Then by using $\eqref{eq:2.16}$ we compute the energy in each of the other cells to be
$\left(\frac{5}{3}\right)^{n_1}(m_1(x)-m_2(x)-1)^2\mathcal{E}(h_0^y)\allowbreak=\left(\frac{5}{3}\right)^{n_1}\left(\frac{3m_0(x)}{2m_0(x)+1}\right)^2\mathcal{E}(h_0^y),$ and by adding we obtain $\eqref{eq:2.22}$. Then $\eqref{eq:2.23}$ follows by Lemma 2.3.

To obtain the estimate $\eqref{eq:2.24}$ we observe that since $0\leq m_0(x)\leq \frac{3}{10}$ it follows from $\eqref{eq:2.7}$ that $m_0(x)$ is bounded above and below by multiples of $\left(\frac{5}{3}\right)^{n_1-n_2}$.  It follows from $\eqref{eq:2.21}$ that $\mathcal{E}(h_0^x)$ is bounded above and below by multiples of $\left(\frac{5}{3}\right)^{n_1}$ since the infinite series is dominated by its first term. We get the same estimate for $\mathcal{E}(h_1^x)$ using $\eqref{eq:2.22}$ since the second summand is bounded by a multiple of $\left(\frac{5}{3}\right)^{2(n_1-n_2)}\left(\frac{5}{3}\right)^{n_1}\left(\frac{5}{3}\right)^{n_2-n_1}$. Then $\eqref{eq:2.24}$ follows from this estimate and $\eqref{eq:2.23}$. $\Box$

\textbf{Corollary 2.5.} \textit{Let $h$ be the harmonic function on $\Omega_x$ with boundary values $h(q_0)=a$ and $h|_{S(x)}=f$, where
\begin{equation}
f=b+\sum_{\omega}c_\omega\psi_\omega\label{eq:2.27}
\end{equation}
for
\begin{equation}
c_\omega=\int_{S(x)}f\psi_{\omega}d\mu.\label{eq:2.28}
\end{equation}
Then $\mathcal{E}(h)$ is bounded above and below by multiplies of
\begin{equation}
\left(\frac{5}{3}\right)^{n_1}(a-b)^2+\sum_{m=0}^{\infty}\sum_{|\omega|=m}2^m\left(\frac{5}{3}\right)^{n_{m+1}}|c_\omega|^2.\label{eq:2.29}
\end{equation}
In particular, $h$ has finite energy if and only if $\eqref{eq:2.29}$ is finite.}

\textit{Proof.} By subtracting a constant we may assume without loss of generality that $a=0$ (this does not change $c_\omega$). Then from $\eqref{eq:2.27}$ we have
\begin{equation}
h=b(1-h_0)+\sum_{m=0}^\infty\sum_{|\omega|=m}c_\omega h_\omega,\label{eq:2.30}
\end{equation}
and the functions $h_0\cup\{h_\omega\}$ are orthogonal in energy by symmetry considerations. Thus
\begin{equation}
\mathcal{E}(h)=b^2\mathcal{E}(1-h_0)+\sum_{m=0}^\infty\sum_{|\omega|=m}|c_\omega|^2\mathcal{E}(h_\omega)\label{eq:2.31}
\end{equation}
and the result follows by the estimates $\eqref{eq:2.24}$. $\Box$

We are also interested in the corresponding result for the $L^2$ norm of $h$. Using similar reasoning we can show that $\parallel h\parallel_2^2$ is bounded above and below by multiples of
\begin{equation}
\left(\frac{1}{3}\right)^{n_1}(a^2+b^2)+\sum_{m=0}^\infty\sum_{|\omega|=m}2^m\left(\frac{1}{3}\right)^{n_{m+1}}|c_\omega|^2.\label{eq:2.32}
\end{equation}
Of course this allows the coefficients to grow so that $\sum_{\omega}|c_\omega|^2$ is infinite, meaning that the boundary values $f$ on $S(x)$ may not be in $L^2(S(x))$.

\section{Normal Deriatives}

We follow the general outline from {[}OS{]} to define normal derivative on $S(x)$. We define
\begin{equation}
\partial_n u|_{S(x)}=\lim_{m\rightarrow\infty}2^m\sum_{|\omega|=m}(-\partial_n(\tilde{F}_\omega q_0))\chi_{S(x)\cap \tilde{F}_\omega(\mathcal{SG})}\label{eq:3.1}
\end{equation}
if the limit exists, where
\begin{equation}
\tilde{F}_\omega=F_0^{n_1-1}F_{\omega_1}F_0^{n_2-n_1-1}F_{\omega_2}\cdots F_0^{n_m-n_{m-1}-1}F_{\omega_m}\label{eq:3.2}.
\end{equation}
The cells $\tilde{F}_\omega(\mathcal{SG})$ for $|\omega|=m$ cover $S(x)$, and $\tilde{F}_\omega q_0$ is the top vertex. Since $\partial_n u(\tilde{F}_\omega q_0)$ is outer directed, upward, we insert the minus sign to get an outer directed normal along $S(x)$.

\textbf{Lemma 3.1.} \textit{$\partial_n h_0^x$ is the constant function on $S(x)$ with value $-2\left(\frac{5}{3}\right)^{n_1}\linebreak[4](1-m_0(x)).$}

\textit{Proof.} We compute $\partial_n h_0^x(q_0)=2\left(\frac{5}{3}\right)^{n_1}(1-m_0(x))$ from the cell $F_0^{n_1}(\mathcal{SG})$. Next consider the cell $F_0^{n_1-1}F_{\omega_1}F_0^{n_2-n_1}(\mathcal{SG})$. The top vertex is $F_{0}^{n_1-1}F_{\omega_1}q_0$, and by symmetry (on the cell $F_0^{n_1}(\mathcal{SG})$), $\partial_n h_0^x(F_0^{n_1-1}F_{\omega_1}q_0)=\frac{1}{2}\partial_n h_0^x(q_0)$ for $\omega_1=1,2$. Thus $2\sum_{|\omega|=1}(-\partial_n h_0^x(\tilde{F}_\omega q_0))\chi_{S(x)\cap \tilde{F}_\omega(\mathcal{SG})}=-\partial_n h_0^x(q_0)\chi_{S(x)}$. By similar reasoning there is no change on the right side of $\eqref{eq:3.1}$ as $m$ increases. $\Box$

\textbf{Lemma 3.2.} \textit{$\partial_n h_{\omega}^x=6\cdot 2^m\left(\frac{5}{3}\right)^{n_{m+1}}\left(\frac{1-m_0(y_m)}{2m_0(y_m)+1}\right)\psi_\omega$.}

\textit{Proof.} On the cell $F_0^{n_1}(\mathcal{SG})$ we compute (using $\eqref{eq:2.14}$)
$$\partial_n h_1^x(F_0^{n_1-1}F_{1}q_0)=-\partial_n h_1^x(F_0^{n_1-1}F_2 q_0)=3\left(\frac{5}{3}\right)^{n_1}(m_1(x)-m_2(x)),$$ so on the cell $F_0^{n_1-1}F_{\omega_1}F_0^{n_2-n_1}(\mathcal{SG})$ we have $$\sum_{|\omega|=1}2(-\partial_n h_1^x(\tilde{F}_\omega q_0))\chi_{S(x)\cap \tilde{F}_\omega(\mathcal{SG})}=6\left(\frac{5}{3}\right)^{n_1}(m_1(x)-m_2(x))\psi_\emptyset,$$
and by the same reasoning as in Lemma 3.1, this does not change as we increase $m$. So this gives the correct result for $\omega=\emptyset$. We then use Lemma 2.3 to scale the result for general $\omega$. $\Box$

\textbf{Theorem 3.3.} \textit{Suppose $h$ and $f$ are given as in Corollary 2.5. Then $\partial_n h$ is given by
\begin{equation}
2(b-a)\left(\frac{5}{3}\right)^{n_1}(1-m_0(x))+\sum_{m=0}^\infty\sum_{|\omega|=m}6\cdot 2^m\left(\frac{5}{3}\right)^{n_{m+1}}\left(\frac{1-m_0(y_m)}{2m_0(y_m)+1}\right)c_\omega\psi_\omega.\label{eq:3.3}
\end{equation} In other words, the Dirichlet-to-Neumann map $f\rightarrow\partial_n h$ is a Haar series multiplier map with multiplier $6\cdot 2^m \left(\frac{5}{3}\right)^{n_{m+1}}\left(\frac{1-m_0(y_m)}{2m_0(y_m)+1}\right)$.}

\textbf{Corollary 3.4.} \textit{Suppose $f$ satisfies
\begin{equation}
\sum_{m=0}^\infty \sum_{|\omega|=m}2^{2m}\left(\frac{5}{3}\right)^{2n_{m+1}}|c_\omega|^2<\infty.\label{eq:3.4}
\end{equation}
Then $\partial_n h$ is well-defined in $L^2(S(x))$ and $\parallel\partial_n h\parallel_2^2$ is bounded above and below by a multiple of $\eqref{eq:3.4}$.}

\textit{Proof.} The Theorem follows from Lemma 3.2, and the Corollary follows from the fact that $\frac{1-m_0(x)}{2m_0(x)+1}$ is uniformly bounded above and below independent of $x$. $\Box$

Note that the finiteness of $\eqref{eq:3.4}$ is a stronger condition than the finiteness of $\eqref{eq:2.29}$, so harmonic functions of finite energy do not necessarily satisfy $\eqref{eq:3.4}$, but functions $h$
satisfying the conditions of Corollary 3.4 automatically have finite energy.

\textbf{Corollary 3.5.} \textit{Suppose $h$ satisfies the hypothesis of Corollary 3.4, and $v$ is any function of finite energy of $\Omega_x$, then the following Gauss-Green formula holds:
\begin{equation}
\mathcal{E}(h,v)=v(q_0)\partial_n h(q_0)+\int_{S(x)}v\partial_n h d\mu\label{eq:3.5}.
\end{equation}}

\textit{Proof.} Apply the standard Gauss-Green formula on the domain $\bigcup_{|\omega|\leq m}\tilde{F}_\omega(\mathcal{SG})$ and take the limit as $m\rightarrow\infty$. $\Box$

\section{Extending Functions of Finite Energy}

In this section we will write $\Omega_x^+$ for the region above $L(x)$ that was previously denoted $\Omega_x$, and $\Omega_x^-$ for the region below $L(x)$.  Under the assumption that $x$ is not a dyadic rational, $S(x)$ is the common boundary of $\Omega_x^+$ and $\Omega_x^-$. The first issue that we address is under what conditions can we glue together functions $u^{\pm}$ of finite energy on $\Omega_x^{\pm}$ to obtain a function of finite energy on $\mathcal{SG}$. Since functions of finite energy are continuous, $u^{\pm}$ must have boundary values on $S(x)$ that agree. It turns out that this is the only condition that we need to impose. This is not surprising since the same is true for gluing functions of finite energy on domains that intersect at a finite set of points.

\textbf{Theorem 4.1.} \textit{Let $u^\pm\in dom\mathcal{E}_{\Omega_x^{\pm}}$, and suppose
\begin{equation}
u^+|_{S(x)}=u^-|_{S(x)},\label{eq:4.1}
\end{equation}
the values being defined by continuity. Then
\begin{eqnarray}\label{eq:4.2}
u=\left\{\aligned
u^+ \mbox{ on } \overline{\Omega}_x^+,\\
u^-\mbox{ on }\overline{\Omega}_x^-,
\endaligned\right.\end{eqnarray}
belongs to $dom\mathcal{E}$ in $\mathcal{SG}$ and
\begin{equation}
\mathcal{E}(u)=\mathcal{E}_{\Omega_x^+}(u^+)+\mathcal{E}_{\Omega_x^-}(u^-).\label{eq:4.3}
\end{equation}}

\textit{Proof.} Let $S_m$ denote the strip of $2^m$ cells of order $n_m$ containing $S(x)$, and let $B_m^\pm$ denote the unions of the cells of order $n_m$ contained in $\Omega_x^\pm$. Then
$$
\mathcal{E}^{(n_m)}(u)=\mathcal{E}_{B_m^+}^{(n_m)}(u)+\mathcal{E}_{B_m^-}^{(n_m)}(u)+\mathcal{E}_{S_m}^{(n_m)}(u).
$$Since $\mathcal{E}_{B_m^\pm}^{(n_m)}(u)\rightarrow\mathcal{E}_{\Omega_x^\pm}(u^\pm)$ as $m\rightarrow\infty$, it suffices to show
\begin{equation}
\mathcal{E}_{S_m}^{(n_m)}(u)\rightarrow 0.\label{eq:4.4}
\end{equation}

Let $C$ denote one of the $n_m$-cells in $S_m$ with boundary points $x_m\in\Omega_x^+$ and $y_m,z_m\in\Omega_x^-$. We need to estimate
\begin{equation}
\left(\frac{5}{3}\right)^{n_m}[(u^+(x_m)-u^-(y_m))^2+(u^+(x_m)-u^-(z_m))^2+(u^-(y_m)-u^-(z_m))^2].\label{eq:4.5}
\end{equation}
It suffices to estimate the first two terms in $\eqref{eq:4.5}$
since $u^-(y_m)-u^-(z_m)=(u^+(x_m)-u^-(z_m))-(u^+(x_m)-u^-(y_m))$,
and by symmetry it suffices to estimate the first term. Let
$S_m^\pm$ be the portion of $S_m$ above or below $S(x)$. There will
be an infinite sequence of points $\{x_m, x_{m+1},\cdots\}$  in
$S_m^+$ and $\{y_m,y_{m+1},\cdots\}$ in $S_m^-$, both converging to
the same point $p\in S(x)$. Since $u^+(p)=u^-(p)$ by
$\eqref{eq:4.1}$, we may write
\begin{equation}
u^+(x_m)-u^-(y_m)=\sum_{j=m}^\infty(u^+(x_j)-u^+(x_{j+1}))-\sum_{j=m}^\infty(u^-(y_j)-u^-(y_{j+1})).\label{eq:4.6}
\end{equation}
Now each pair $(x_j,x_{j+1})$ are vertices of a cell $C_j$ of order $n_{j+1}$ in $\Omega_x^+$. Note that all these cells are essentially disjoint.

\begin{figure}[h]
\begin{center}
\includegraphics[width=10cm,totalheight=8.7cm]{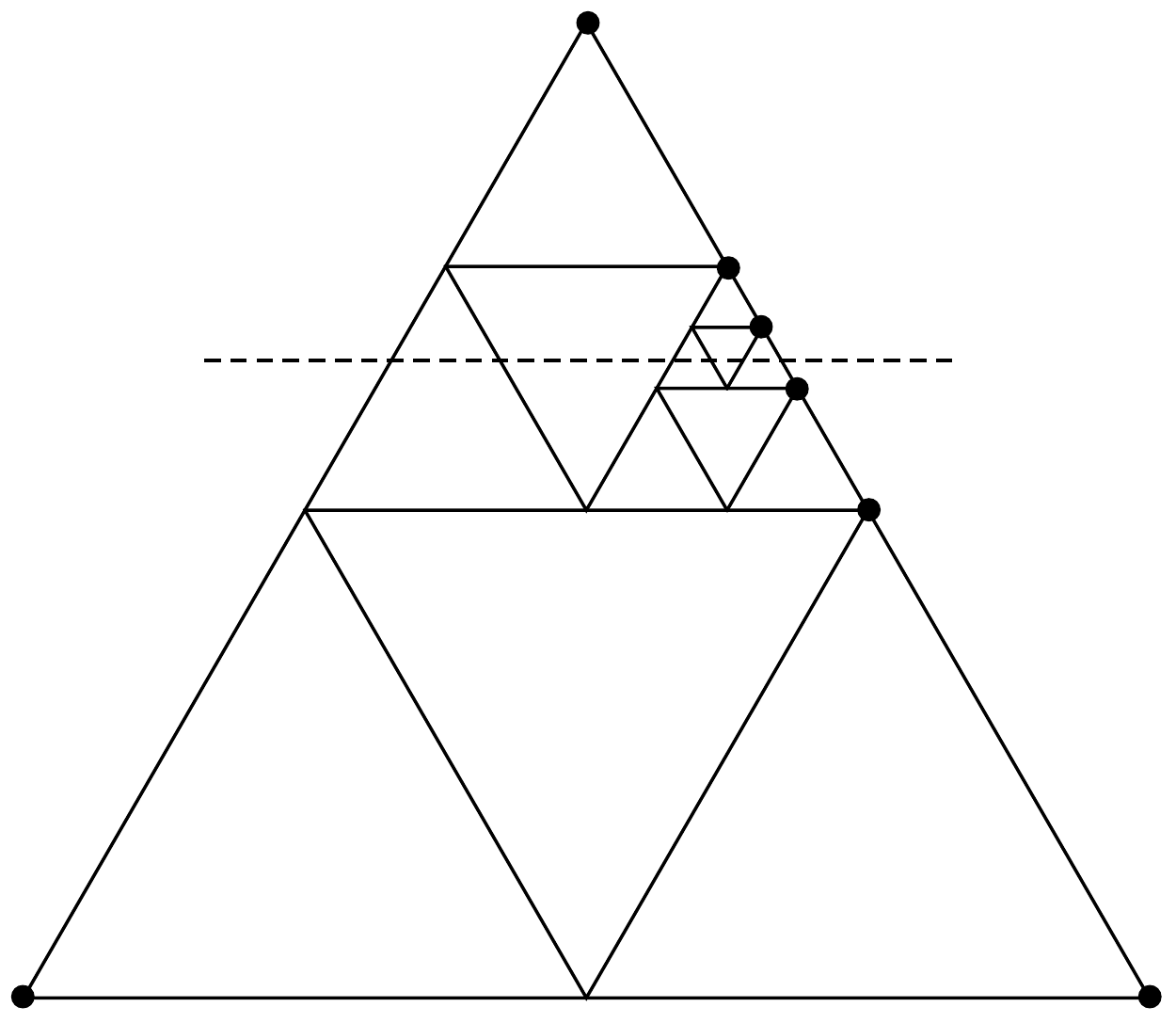}
\setlength{\unitlength}{1cm}
\begin{picture}(0,0) \thicklines
\put(-2.6,4.3){$y_{m+1}$}
\put(-3.2,5.3){$y_{m+2}$}
\put(-3.6,5.9){$x_{m+2}$}
\put(-3.8,6.4){$x_{m+1}$}
\put(-5.4,8.8){$x_m$}
\put(-.2,0.2){$y_m$}
\put(-10.7,0.2){$z_m$}
\put(-2,5.5){$S(x)$}
\end{picture}
\begin{center}
\textbf{Figure 4.1.}
\end{center}
\end{center}
\end{figure}

So we have the estimate
\begin{equation}
|u^+(x_j)-u^+(x_{j+1})|\leq\left(\frac{3}{5}\right)^{n_{j+1}/2}\mathcal{E}_{C_j}(u^+)^{1/2}. \label{eq:4.7}
\end{equation}
By the Cauchy-Schwarz inequality we obtain
\begin{eqnarray}
\sum_{j=m}^\infty|u^+(x_j)-u^+(x_{j+1})|&\leq&\left(\sum_{j=m}^\infty\left(\frac{3}{5}\right)^{n_{j+1}}\right)^{1/2}\left(\sum_{j=m}^\infty\mathcal{E}_{C_j}(u^+)\right)^{1/2}\label{eq:4.8}\\&\leq& c\left(\frac{3}{5}\right)^{n_{m}/2}\mathcal{E}_{C\cap S_m^+}(u^+)^{1/2}.\nonumber
\end{eqnarray}
By similar reasoning we obtain the same estimate with $|u^-(y_j)-u^-(y_{j+1})|$, so by $\eqref{eq:4.6}$ we have
\begin{equation}
\left(\frac{5}{3}\right)^{n_m}|u^+(x_m)-u^-(y_m)|^2\leq c\mathcal{E}_{C\cap S_{m}^+}(u^+)+c\mathcal{E}_{C\cap S_m^-}(u^-).\label{eq:4.9}
\end{equation}
Summing $\eqref{eq:4.9}$ over all the $2^m$ cells $C$ yields
\begin{equation}
\mathcal{E}_{S_m}^{(n_m)}(u)\leq c\mathcal{E}_{S_m^+}(u^+)+c\mathcal{E}_{S_m^-}(u^-)\label{eq:4.10}
\end{equation}
and $\mathcal{E}_{S_m^\pm}(u^\pm)\rightarrow 0$ because $\bigcap_m S_m^\pm=S(x)$ and $S(x)$ has measure zero in the Kusuoka measure. $\Box$

It is easy to characterize the restrictions to $S(x)$ of functions of finite energy on $\Omega_x^+$.

\textbf{ Theorem 4.2.} \textit{A function $f$ on $S(x)$ is the restriction to $S(x)$ of a function $u^+$ of finite energy on $\Omega_x^+$ if and only if $f$ has a Haar series expansion $\eqref{eq:2.27}$
 with $\eqref{eq:2.29}$ finite (here a=0), and $\eqref{eq:2.29}$ is bounded by a multiple of $\mathcal{E}_{\Omega_x^+}(u^+)$.}

\textit{Proof.} Let $h$ be the harmonic function on $\Omega_x^+$ with the same boundary values $f$. Since harmonic functions minimize energy, $\mathcal{E}_{\Omega_x^+}(h)\leq\mathcal{E}_{\Omega_x^+}(u^+)$, and the result follows from Corollary 2.5. $\Box$

However, there is no such simple result for $\Omega_x^-$. We pose the following extension problem.

\textbf{Problem 4.3.} \textit{Does there exist a bounded linear extension operator (meaning $Tu|_{\Omega_x^+}=u$) $T: dom_{\Omega_x^+}(\mathcal{E})\rightarrow dom_{\mathcal{SG}}(\mathcal{E})$?}

There is a simple obstruction to solving this problem.

\textbf{Definition 4.4.} $x$ satisfies the \textit{nonconsecutive condition} with bound $N$ if there are no $N$ consecutive integers in the sequence $\{n_m\}$.

Note that a generic value of $x$ will not satisfy this condition.
However there are uncountably many (of Hausdorff dimension $1$)
values of $x$  that do satisfy the condition. Perhaps the simplest
choice has $n_m=2m-1$, with $N=2$.

\textbf{Theorem 4.5.} \textit{ Let $E$ denote the collection of
$x$ satisfying the nonconsecutive condition. Then the Hausdorff
dimension of $E$ is $1$.}

\textit{Proof.} Let $E_N$ denote the collection of $x$ satisfying
the nonconsecutive condition with bound $N$. Then $E=\bigcup_{N\geq
2} E_N$ and
\begin{equation}
E_2\subset E_3\subset \cdots\subset E_N\subset\cdots.\label{eq:4.11}
\end{equation}
We will first prove that the
Hausdorff dimension of $E_N$ is the unique positive root of the
equation
\begin{equation}
2-2^s-2^{-Ns}=0.\label{eq:4.12}
\end{equation}

Consider the set $E_N$. We divide it into the disjoint union
$E_N=\bigcup_{k\geq 1}E_{N,k}$ where $E_{N,k}$ is the set of $x$ in
$E_N$ whose $n_1$-digit is $k$. Obviously, for each $k$, $E_{N,k}$
is a similar copy of $E_{N,1}$ with contraction ratio $2^{1-k}$.
Since the Hausdorff dimension is stable for countable unions, we
just need to compute the dimension of $E_{N,1}$. For this set, by
the nonconsecutive condition, we can write
\begin{equation}
E_{N,1}=\left(\bigcup_{j\geq
3}(2^{-1}+E_{N,j})\right)\cup\cdots\cup\left(\bigcup_{j\geq{N+1}}(2^{-1}+\cdots+2^{-(N-1)}+E_{N,j})\right).\label{eq:4.13}
\end{equation}
Since $|E_{N,j}|\leq 1/2^j$, it is easy to check that the above
union is disjoint. Moreover, $\eqref{eq:4.13}$ is essentially a
self-similar identity for the set $E_{N,1}$ with contraction ratios,
$$2^{-2},2^{-3},\cdots; 2^{-3},2^{-4},\cdots;2^{-N},2^{-(N+1)},\cdots,$$
satisfying the open set condition ( with the open set $(2^{-1},1)$).
(See {[}M{] for the theory of infinitely generated self-similar
sets.) Hence the Hausdorff dimension of $E_{N,1}$ is the solution of
the equation
\begin{equation}
1=\sum_{k=2}^N\sum_{j\geq k}(2^{-s})^j=\sum_{k=2}^N\frac{(2^{-s})^k}{1-2^{-s}}=\frac{2^{-2s}-2^{-s(N+1)}}{(1-2^{-s})^2},\label{eq:4.14}
\end{equation}
which simplifies  to $\eqref{eq:4.12}$. So we get the
Hausdorff dimension of $E_N$.

Using $\eqref{eq:4.11}$, an easy calculation will show that the
Hausdorff dimension of $E$ is $1. \Box$

If $x$ fails to satisfy this condition, then there are pairs of points in $\Omega_x^+$ that are much closer to each other in $\mathcal{SG}$ than in $\Omega_x^+$. For example, if $n_j=j$ for $j\leq N$ then the points $F_1F_2^{N-1}q_0$ and $F_2F_1^{N-1}q_0$ in $\Omega_x^+$ are distance on the order of $\left(\frac{3}{5}\right)^N$ apart in the resistance metric on $\mathcal{SG}$, but are far apart in $\Omega_x^+$. Note that $h_1^x(F_1F_2^{N-1}q_0)-h_1^x(F_2F_1^{N-1}q_0)=2h_1^x(F_1F_0^{N-1}q_0)$ and $\mathcal{E}(h_1^x)$ is bounded. The estimate analogous to $\eqref{eq:4.7}$ shows
$$c\leq\left(\frac{3}{5}\right)^{N/2}\mathcal{E}(u)^{1/2}$$
for any extension $u$ of $h_1^x$ to $\mathcal{SG}$, hence $\mathcal{E}(u)\geq c\left(\frac{5}{3}\right)^N$. This means that the bound on the operator $T$, if it exists, would be bounded below by a multiple of $\left(\frac{5}{3}\right)^{N/2}$.

\begin{figure}[h]
\begin{center}
\includegraphics[width=6cm,totalheight=5.3cm]{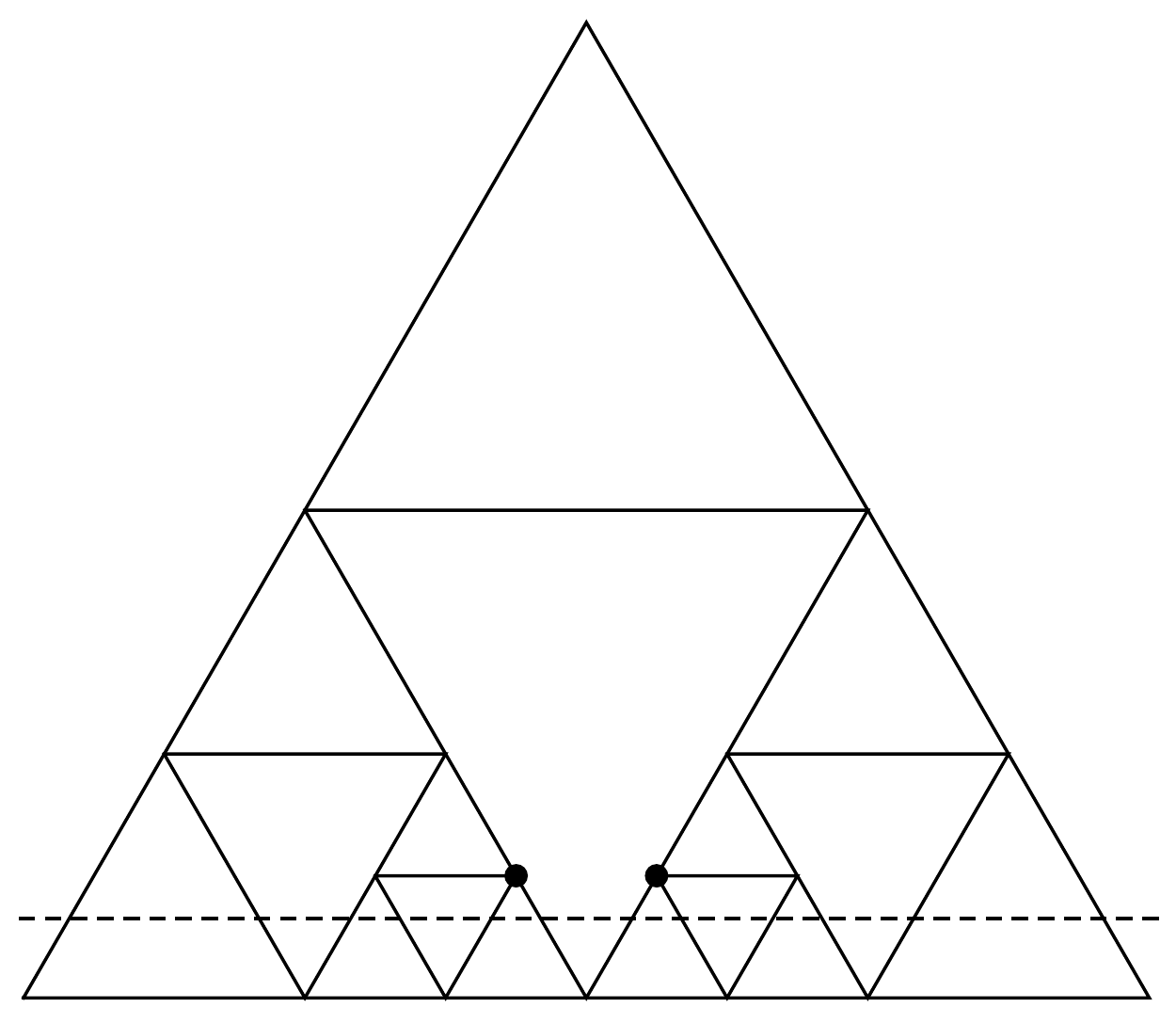}
\setlength{\unitlength}{1cm}
\begin{picture}(0,0) \thicklines
\put(-0.1,0.5){$S(x)$}
\end{picture}
\begin{center}
\textbf{Figure 4.2.}
\end{center}
\end{center}
\end{figure}

The same reasoning applies locally if $\{n_m\}$ has a consecutive string of $N$ integers. Thus if such strings exist for all $N$ then $T$ cannot be bounded.  On the other hand it is easy to see that if the nonconsecutive condition holds for $x$ then distances in $\Omega_x^+$ and $\mathcal{SG}$ are comparable. Note that this is very reminiscent of the type of condition that appears in the work of Peter Jones in the Sobolev extension problem in domains in Euclidean space ({[}J{]}, {[}R{]}).

\textbf{Theorem 4.6.} \textit{The extension Problem 4.3 has a
positive solution if and only if $x$ satisfies the nonconsecutive
condition, in which case the bound on $T$ is
$O(\left(\frac{10}{3}\right)^{N/2})$.}

\textit{Proof.} We need to construct an extension operator $T$ under the assumption that $x$ satisfies the nonconsecutive condition. In view of our previous results, it suffices to solve the extension problem for the functions $h_\omega$ (and also $1-h_0$), say $T{h_\omega}=\tilde{h}_\omega$ where the function $\tilde{h}_\omega$ are orthogonal in energy and
\begin{equation}
\mathcal{E}(\tilde{h}_\omega)\leq
C(N)2^m\left(\frac{5}{3}\right)^{n_{m+1}}.\label{eq:4.15}
\end{equation}

Suppose first that $N=2$. Consider first $1-h_0$ and $h_\emptyset$. Assume for simplicity that $n_1=1$. Then $n_2\geq 3$. Then $S(x)$ passes through the cells $F_1F_0(\mathcal{SG})$ and $F_2F_0(\mathcal{SG})$. We will extend $1-h_0$ to be identically $1$ on the bottom portions of theses cells, and $h_\emptyset$ to be $1$ on $F_1F_0(\mathcal{SG})$ and $-1$ on $F_2F_0(\mathcal{SG})$. On the remaining four cells of level $2$ we make the extension harmonic with boundary values $0$ on the bottom vertices (see Figure 4.3).

\begin{figure}[h]
\begin{center}
\includegraphics[width=11.8cm,totalheight=4.2cm]{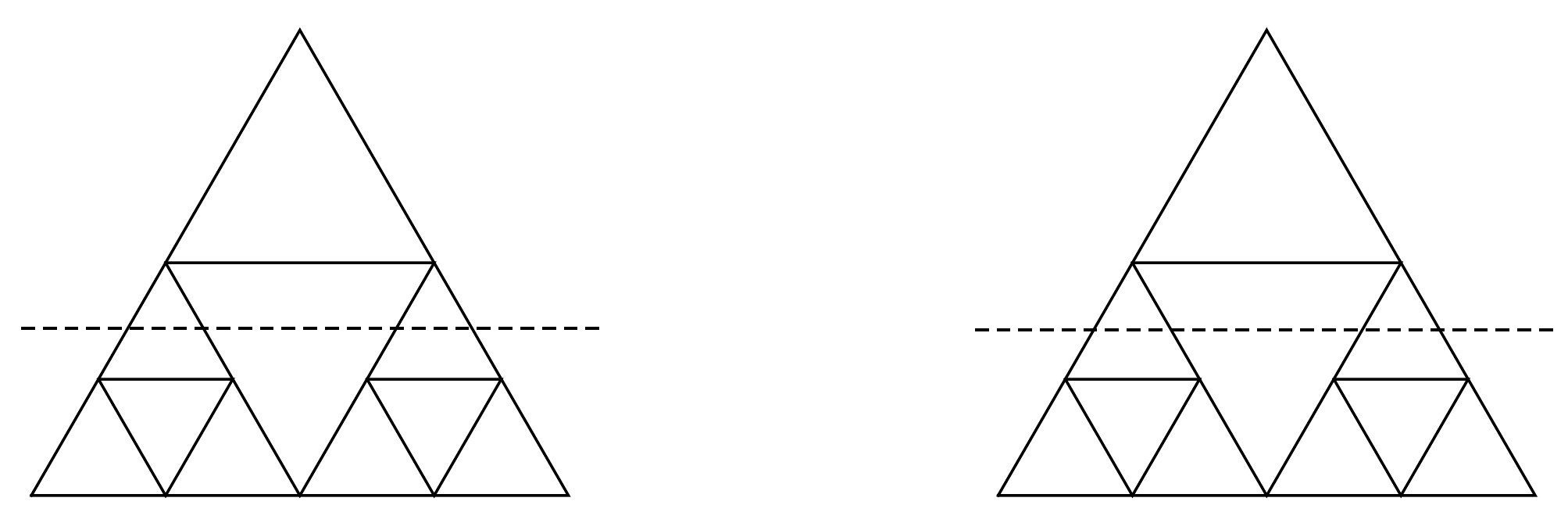}
\setlength{\unitlength}{1cm}
\begin{picture}(0,0) \thicklines
\put(-0.3,1.5){$S(x)$}\put(-11.9,0){$0$}\put(-10.9,0){$0$}
\put(-7.5,1.5){$S(x)$}\put(-9.85,0){$0$}
\put(-10.9,1.3){$1$}\put(-8.85,0){$0$}
\put(-8.85,1.3){$1$}\put(-7.85,0){$0$}
\put(-1.7,1.3){$-1$}\put(-4.6,0){$0$}
\put(-3.6,1.3){$1$}\put(-3.6,0){$0$}
\put(-0.9,1.1){$-1$}\put(-2.6,0){$0$}
\put(-2.5,1.1){$-1$}\put(-1.6,0){$0$}
\put(-3.0,1.1){$1$}\put(-0.55,0){$0$}
\put(-4.2,1.1){$1$}\put(-2.6,4){$0$}
\put(-8.2,1.1){$1$}\put(-9.85,4){$0$}
\put(-9.5,1.1){$1$}
\put(-10.25,1.1){$1$}
\put(-11.5,1.1){$1$}
\put(-11,-0.5){extension of $1-h_0$}
\put(-3.5,-0.5){extension of $h_{\emptyset}$}
\end{picture}\vspace{1cm}
\begin{center}
\textbf{Figure 4.3.}
\end{center}
\end{center}
\end{figure}

Note that the added energy of these extensions is exactly $8\left(\frac{5}{3}\right)^2$. Also, since
 one extension is symmetric and one is skew-symmetric with respect to the vertical reflection, they are orthogonal in energy. If $n_1>1$ we may repeat the same process on $F_0^{n_1-1}(\mathcal{SG})$ and then continue the extension to be identically zero on the complement of $F_0^{n_1-1}(\mathcal{SG})$. The added energy is exactly $8\left(\frac{5}{3}\right)^{n_1+1}$, but the energy of the original functions was also a multiple of $\left(\frac{5}{3}\right)^{n_1}$, so this is consistent with $\eqref{eq:4.15}$ with $m=0$ and gives a uniform bound on the extension operator.

For the extension of $h_\omega$ we just have to repeat the same procedure miniaturized. If $|\omega|=m$ then $h_\omega$ is supported on a cell of order $n_{m+1}-1$ and since $n_{m+2}\geq n_{m+1}+2$ the right side of Figure 4.3 describes $h_\omega$ and its extension (except for a factor of $2^{m/2}$) to that cell, and then we may glue this to zero in the complement of the cell. Thus we get an extension with the same energy bound. For words $\omega$ with $|\omega|=m$, the extended function have disjoint support, so the energies are orthogonal. Comparing extensions for words of different length with overlapping support, we again have a symmetry/skew-symmetry dichotomy  with respect to the local reflection in the vertical axis of the smaller supporting cell (this is the overlap of the supports) and so we again have energy orthogonality. This completes the proof for $N=2$.

For general $N$ the argument is simlar. In Figure 4.4 we show the extension of $h_\emptyset$ when $N=3$ and $n_1=1, n_2=2, n_3\geq 4$.

\begin{figure}[h]
\begin{center}\vspace{0.3cm}
\includegraphics[width=10cm,totalheight=8.7cm]{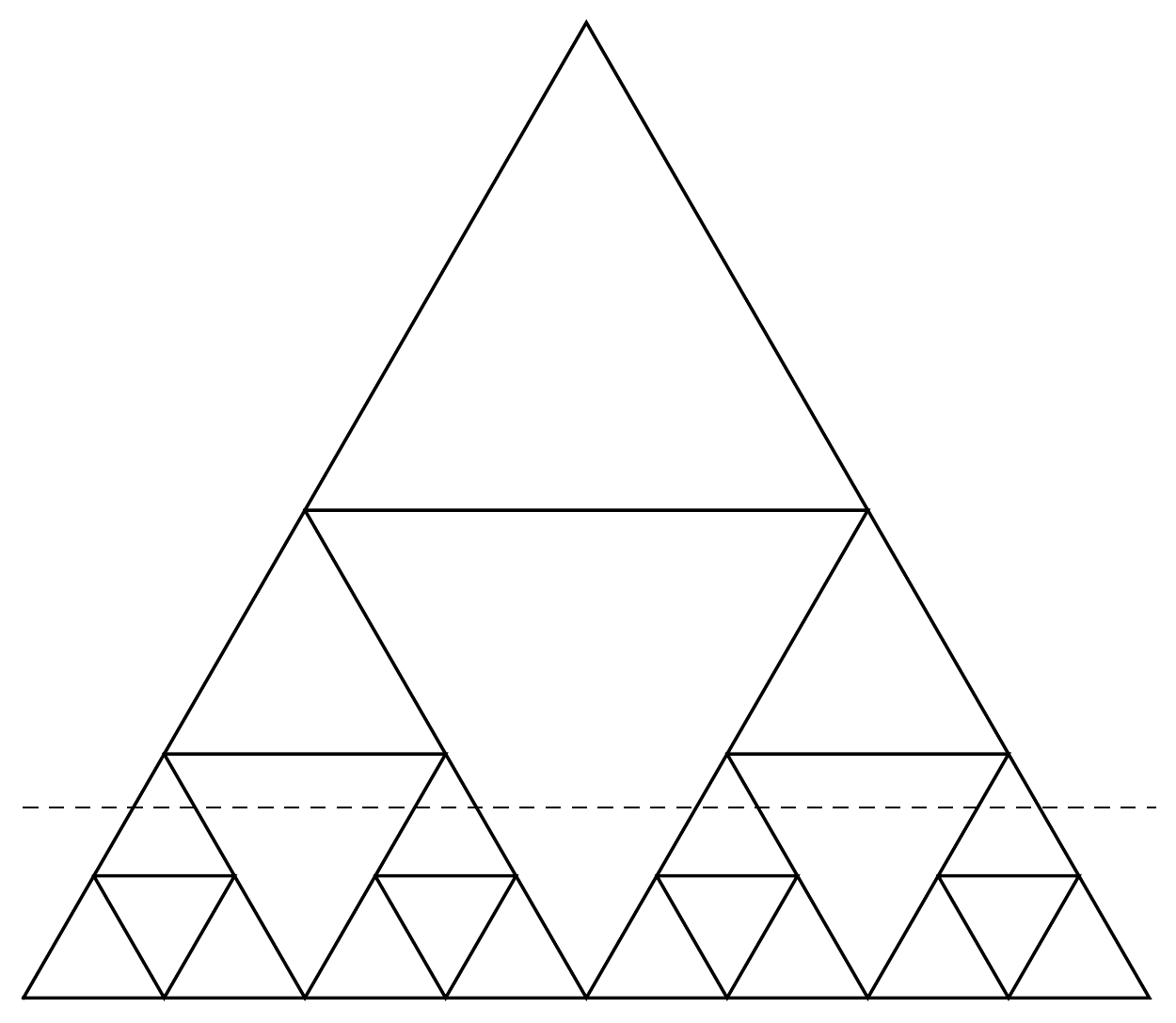}
\setlength{\unitlength}{1cm}
\begin{picture}(0,0) \thicklines
\put(-0.3,1.7){$S(x)$}
\put(-0.5,-0.1){$0$}\put(-1.7,-0.1){$0$}\put(-2.9,-0.1){$0$}
\put(-4.1,-0.1){$0$}
\put(-5.3,-0.1){$0$}
\put(-6.5,-0.1){$0$}
\put(-7.7,-0.1){$0$}\put(-8.9,-0.1){$0$}\put(-10.1,-0.1){$0$}
\put(-0.9,1.2){$-1$}\put(-2.8,1.2){$-1$}\put(-3.4,1.2){$-1$}
\put(-5.1,1.2){$-1$}
\put(-5.8,1.2){$1$}\put(-7.3,1.2){$1$}\put(-8.2,1.2){$1$}\put(-9.7,1.2){$1$}
\put(-8.9,1.5){$1$}\put(-6.5,1.5){$1$}\put(-1.9,1.5){$-1$}\put(-4.3,1.5){$-1$}
\put(-5.3,8.7){$0$}
\end{picture}\vspace{1cm}
\begin{center}
\textbf{Figure 4.4.}
\end{center}
\end{center}
\end{figure}

Here we have $2^N$ cells of order $N$ contributing to the energy, and this multiplies the energy by $O(\left(\frac{10}{3}\right)^N)$. Since the norm of the extension is measured in terms of the square root of the energy, we obtain the $O(\left(\frac{10}{3}\right)^{N/2})$ bound. $\Box$

The optimal extension operator would produce functions that are harmonic on $\Omega_x^-$. In particular, it would be interesting to have an explicit description of the functions $h_\omega^-$ that are harmonic on $\Omega_x^-$ and are equal to $\psi_\omega$ on $S(x)$, again under the nonconsecutive condition.

We may regard Theorem 4.2 as a trace theorem and Theorem 4.6 as an
extension theorem for $dom\mathcal{E}$ regarded as a Sobolev space,
and then we should ask if there are analogous results for other
Sobolev spaces. In {[}S2{]} the spaces $dom_{L^2}(\Delta^k)$ on
$\mathcal{SG}$ are considered as Sobolev spaces
($dom_{L^2}(\Delta^k)=\{u\in L^2(\mathcal{SG}): \Delta^ju\in
L^2(\mathcal{SG})$ for all $j\leq k$\}). Similarly for the space
$\{u\in dom_{L^2}(\Delta^k):\mathcal{E}(\Delta^k u)<\infty\}$. These
spaces are easily characterized in terms of expansions in
eigenfunctions of the Laplacian. A complete theory of the eigenspaces of the Laplacian on $\Omega_1$ is given in  {[}Q{]}.

\textbf{Problem 4.7.} \textit{For each of these Sobolev spaces,
characterize the space of traces on $S(x)$ and restrictions to
$\Omega_x^+$, for $x$ satisfying the nonconsecutive condition.}

It seems plausible that the trace problem may have a solution with a condition similar to $\eqref{eq:2.29}$ for the Haar expansion $\eqref{eq:2.27}$ with different multiples of $|c_\omega|^2$ depending on the Sobolev space.  The restriction problem is likely to be more challenging. It is clear that restrictions of functions in $dom_{L^2}(\Delta^k)$ must satisfy $\Delta^j u\in L^2(\Omega_x^+)$ for $j\leq k$, but that is not sufficient because all harmonic functions automatically have $\Delta^j u=0$. It would seem that the characterization of restriction Sobolev spaces would also have to involve conditions on traces on $S(x)$. Related problems are discussed in  {[}LS{]} and  {[}LRSU{]}.

\section{Green's Function}

For a given $k$, let $V_k$ denote the set of vertices on the
$k$-level graph approximation of $\mathcal{SG}$. For a point $z\in
V_k\setminus V_0$, let $\phi_{z}^k$ denote the piecewise harmonic
spline of level $k$ satisfying $\phi_z^k(t)=\delta_{zt}$ for $t\in
V_k$ and extended harmonically on $\mathcal{SG}$. Notice that
$\phi_{z}^k\in dom_0\mathcal{E}$ because $z\notin V_0$ and it is
supported in the two $k$-cells meeting at $z$. Recall that in the
standard theory (see the books {[}Ki{]} and {[}S3{]}), the Green's
function $G(s,t)$ to solve the Dirichlet problem $-\Delta u=F$ on
$\mathcal{SG}$, subject to the boundary condition $u|_{V_0}=0$ via
an integral transform $\int_{\mathcal{SG}}G(s,t)F(t)dt$, has the
following explicit formula,

\begin{equation}
G(s,t)=\lim_{M\rightarrow\infty}G^M(s,t) \mbox{ (uniform
limit)}\label{eq:5.1}
\end{equation}
with
\begin{equation}
G^M(s,t)=\sum_{k=1}^M\sum_{z,z'\in V_{k}\setminus
V_{k-1}}g(z,z')\phi_z^k(s)\phi_{z'}^k(t),\label{eq:5.2}
\end{equation}
where
\begin{equation}\label{eq:5.3}
g(z,z')=\left\{\begin{array}{ll}
                \frac{3}{10}\left(\frac{3}{5}\right)^k &  \mbox{ for } z=z'\in V_k\setminus V_{k-1}, \\
                \frac{1}{10}\left(\frac{3}{5}\right)^k & \mbox{ for }
                z\neq z'\in V_{k}\setminus V_{k-1}, \mbox{ contained in the same
                } (k-1)\mbox{-cell},\\
                0, & \mbox{ otherwise.}
              \end{array}
\right.
\end{equation}

To get an analogous Green's function on $\Omega_x$, we should first
modify the definition of those piecewise harmonic splines $\phi_z^k$
whose support intersects the boundary $S(x)$ of the domain
$\Omega_x$. More specially, let $\omega$ be a word of symbols
$\{1,2\}$ with $|\omega|=m$ and $z=\tilde{F}_{\omega}(q_0)$. We
redefine $\phi_z^{n_m}$ to be the piecewise harmonic spline with
value $1$ on $z$, $0$ on $V_{n_m}\cap\Omega_x$ and $S(x)$, and
extended harmonically on $\Omega_x$. Obviously the support of
$\phi_z^{n_m}$ is contained in two $n_m$-cells meeting at $z$, with
$\phi_z^{n_m}=h_0^{y_m}\circ \tilde{F}_\omega^{-1}$ on the cell
$\tilde{F}_\omega(\mathcal{SG})$ and with values unchanged on the
other cell, denoted by $\tilde{\tilde{F}}_\omega(\mathcal{SG})$,
where

\begin{equation}\label{eq:5.4}
\tilde{\tilde{F}}_\omega=\left\{\begin{array}{ll}
                F_0^{n_1-1}F_{\omega_1}F_0^{n_2-n_1-1}F_{\omega_2}\cdots F_0^{n_{m-1}-n_{m-2}-1}F_{\omega_{m-1}}F_0^{n_m-n_{m-1}} & \mbox{ for } m\geq 2, \\
                F_0^{n_1} & \mbox{ for } m=1.
              \end{array}
\right.
\end{equation}

\begin{figure}[h]
\begin{center}\vspace{0.3cm}
\includegraphics[width=4cm,totalheight=5.7cm]{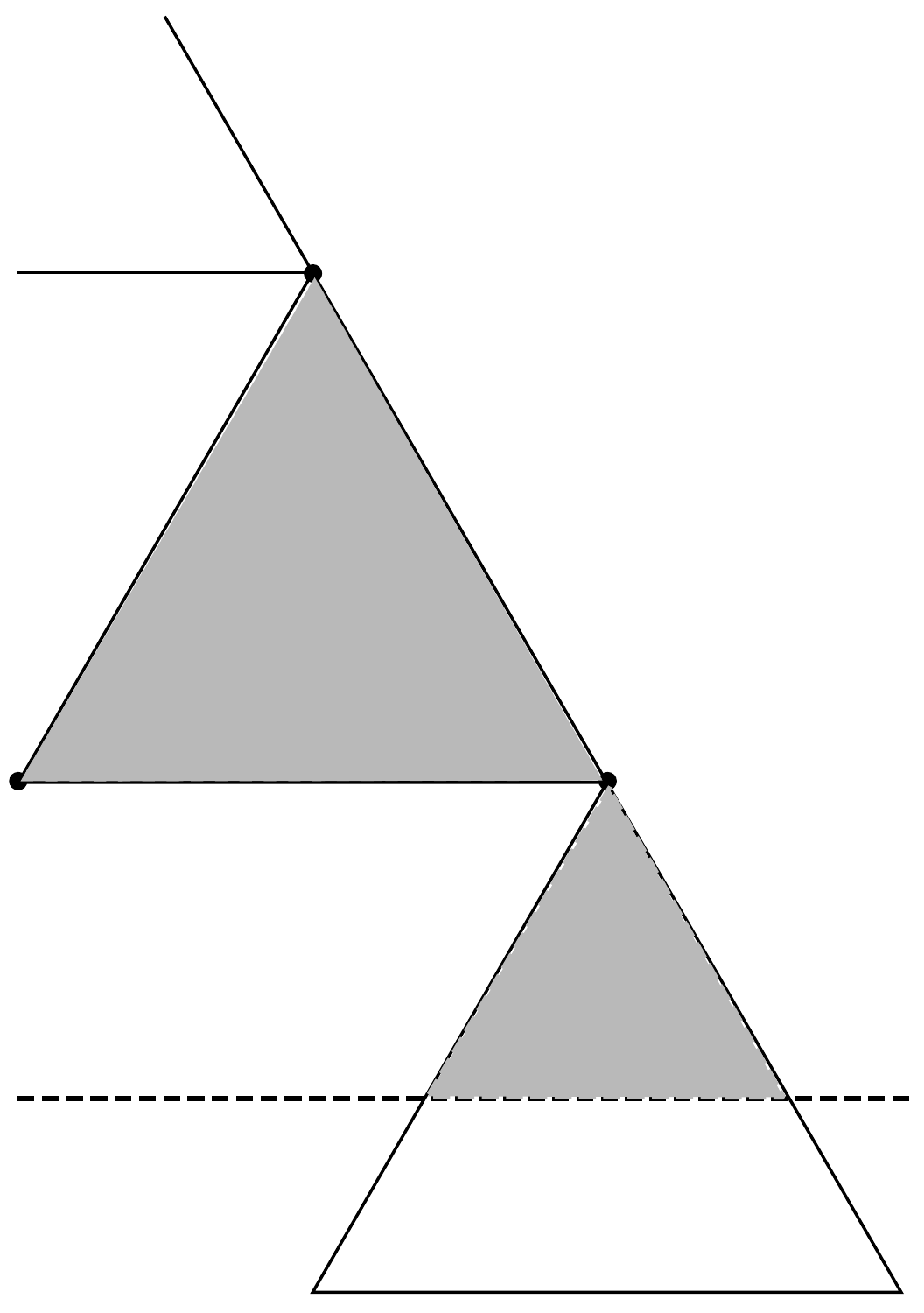}
\setlength{\unitlength}{1cm}
\begin{picture}(0,0) \thicklines
\put(-0.3,0.8){$S(x)$}
\put(-1.5,2.3){$z$}\put(-2.8,4.5){$w$}\put(-4.4,2.3){$z'$}
\put(-2.2,0.8){$\tilde{F}_\omega(\mathcal{SG})$}
\put(-3.4,3){$\tilde{\tilde{F}}_\omega(\mathcal{SG})$}\put(-2,3.3){$n_m$-cell}
\end{picture}
\begin{center}
\textbf{Figure 5.1.}\small{ The support of $\phi_z^{n_m}$.}
\end{center}
\end{center}
\end{figure}

\textbf{Lemma 5.1.} \textit{Let $z=\tilde{F}_{\omega}(q_0)$, then
\begin{equation}\label{eq:5.5}
\mathcal{E}_{\Omega_x}(\phi_z^{n_m},v)=\left(\frac{5}{3}\right)^{n_m}\left(\frac{1+m_0(y_{m-1})}{m_0(y_{m-1})}v(z)-v(z')-v(w)\right)
\end{equation}
for any $v\in dom_0\mathcal{E}_{\Omega_x}$, where
$z'=\tilde{F}_{\omega_1\cdots\omega_{m-1}(3-\omega_m)}(q_0)$ and
$w=\tilde{F}_{\omega_1\cdots\omega_{m-1}}(q_0)$ are the two
$n_m$-neighbors of $z$(See Figure 5.1.).}

\textit{Proof.} On the cell $\tilde{F}_\omega(\mathcal{SG})$, by
using the localized Gauss-Green formula (see $\eqref{eq:3.5}$),
\begin{equation}\mathcal{E}_{\Omega_x\cap
\tilde{F}_\omega(\mathcal{SG})}(\phi_z^{n_m},v)=v(z)\partial_n\phi_z^{n_m}(z)=2\left(\frac{5}{3}\right)^{n_{m+1}}(1-m_0(y_m))v(z).\label{eq:5.6}
\end{equation}
The last equality follows from the same argument as the proof of
Lemma 3.1 with suitable scaling.

On the other cell $\tilde{\tilde{F}}_\omega(\mathcal{SG})$, by using
the standard theory, \begin{equation}\mathcal{E}_{\Omega_x\cap
\tilde{\tilde{F}}_\omega(\mathcal{SG})}(\phi_z^{n_m},v)=\left(\frac{5}{3}\right)^{n_{m}}(2v(z)-v(z')-v(w)).\label{eq:5.7}\end{equation}

Summing the energies on the two cells, we get the desired result by
using $\eqref{eq:2.5}. \Box$

Let $T_x^m$ be the set of vertices in $V_{n_m}\cap\Omega_x$ which
can be expressed as $\tilde{F}_\omega(q_0)$ for some word
$\omega=\omega_1,\cdots,\omega_m$ of symbols $\{1,2\}$, and
$T_x=\bigcup_{m\geq 1}T_x^m$.

\textbf{Definition 5.2.} For fixed $m$, let
\begin{equation}\label{eq:5.8}
G_{\Omega_x}^m(s,t)=\sum_{k=1}^{n_m}\sum_{z,z'\in (V_{k}\setminus
V_{k-1})\cap \Omega_x}g_x(z,z')\phi_z^k(s)\phi_{z'}^k(t),
\end{equation}
with
\begin{equation}\label{eq:5.9}
g_x(z,z')=\left\{\begin{array}{ll}
                \frac{m_0(y_{l-1})+m_0(y_{l-1})^2}{2m_0(y_{l-1})+1}\left(\frac{3}{5}\right)^{n_{l}} &  \mbox{ for } z=z'\in T_x^l \mbox{ with } l\leq m, \\
                \frac{m_0(y_{l-1})^2}{2m_0(y_{l-1})+1}\left(\frac{3}{5}\right)^{n_{l}} & \mbox{ for }
                z\neq z'\in T_x^l, \mbox{ being } n_l\mbox{-neighbors, with } l\leq m, \\
                g(z,z') & \mbox{ for } z,z'\in V_k\setminus V_{k-1} \mbox{ contained in }  \\ & \mbox{ a
                } (k-1)\mbox{-cell in } \Omega_x,\\
                0, & \mbox{ otherwise.}
              \end{array}
\right.
\end{equation}

Then it is obvious that $G_{\Omega_x}^m(s,t)$ converges uniformly to
a function $G_{\Omega_x}(s,t)$ as $m$ goes to infinity.

\textbf{Theorem 5.3.} \textit{ $G_{\Omega_x}$ is the Green's
function for $\Omega_x$, namely
\begin{equation}u(s)=\int_{\Omega_x}G_{\Omega_x}(s,t)F(t)dt\label{eq:5.10}
\end{equation}
solves the Dirichlet problem $-\Delta u=F$ on $\Omega_x$ with
$u|_{\partial\Omega_x}=0$, for any continuous F.}

\textit{Proof.} Similar to the $\mathcal{SG}$ case, suppose we could prove \begin{equation}\label{eq:5.11}
\mathcal{E}_{\Omega_x}(G_{\Omega_x}^m(\cdot,t),v)=\sum_{z\in
V_{n_m}\cap \Omega_x}v(z)\phi_z^{n_m}(t)
\end{equation}
for any $v\in dom_{0}\mathcal{E}_{\Omega_x}$.

Then just multiply
$\eqref{eq:5.11}$ by $F(t)$ and integrate, using the standard
arguments to interchange the energy and integral, to obtain
\begin{equation}
\mathcal{E}_{\Omega_x}(u_m,v)=\int_{\Omega_x}F(t)\sum_{z\in
V_{n_m}\cap \Omega_x}v(z)\phi_z^{n_m}(t)dt\label{eq:5.12}
\end{equation}
for
\begin{equation}
u_m(s)=\int_{\Omega_x}G_{\Omega_x}^m (s,t)F(t)dt. \label{eq:5.13}
\end{equation}
Since
\begin{equation}
\sum_{z\in V_{n_m}\cap \Omega_x}v(z)\phi_z^{n_m}(t)\rightarrow
v(t)\label{eq:5.14}
\end{equation}
uniformly as $m\rightarrow\infty$, the right side of
$\eqref{eq:5.12}$ converges to $\int_{\Omega_x}F(t)v(t)dt$, and the
left side converges to $\mathcal{E}_{\Omega_x}(u,v)$ as $m$ goes to
$\infty$. Thus we have
\begin{equation}
\mathcal{E}_{\Omega_x}(u,v)=\int_{\Omega_x}Fvdt \label{eq:5.15}
\end{equation}
for any $v\in dom_0\mathcal{E}_{\Omega_x}$, which yields that
$-\Delta u=F$ with $u|_{\partial\Omega_x}=0$.

 Hence our goal is to prove $\eqref{eq:5.11}$. The
function $G_{\Omega_x}^m(s,t)$, which we regard as a function of
the single variable $s$, could be viewed as a linear combination of
terms $\phi_{z}^k(s)$. Then it is clear that
$\mathcal{E}_{\Omega_x}(G_{\Omega_x}^m(\cdot,t),v)$ is a linear
combination of $v(z)$ for $z\in V_{n_m}\cap\Omega_x$. So we need to
compute the combination coefficient of $v(z)$ for each $z$.

 Let $z_0\in V_{n_m}\cap \Omega_x$. If $z_0\notin T_x$, it is
easy to observe that there exists a cell containing $z_0$ as an
interior point. There are many terms in $G_{\Omega_x}^m$ contribute
to the coefficient of $v(z_0)$, but these terms all have supports
away from $S(x)$. Thus the standard argument for the $\mathcal{SG}$
case shows that the coefficient of $v(z_0)$ should be
$\phi_{z_0}^{n_m}(t)$.

Hence we only need to consider the case that $z_0\in T_x$. We first
do this when $z_0\in T_{x}^m$. Write $z'_0$ the unique
$n_m$-neighbors of $z_0$ in the same level. Then the
only terms in $G_{\Omega_x}^m$ that contribute to the coefficient of
$v(z_0)$ are
$$g_x(z_0,z_0)\phi_{z_0}^{n_m}(s)\phi_{z_0}^{n_m}(t),
 g_x(z_0,z'_0)\phi_{z_0}^{n_m}(s)\phi_{z'_0}^{n_m}(t),$$
$$g_x(z'_0,z_0)\phi_{z'_0}^{n_m}(s)\phi_{z_0}^{n_m}(t), g_x(z'_0,z'_0)\phi_{z'_0}^{n_m}(s)\phi_{z'_0}^{n_m}(t).$$ By Lemma 5.1, the total contribution is
\begin{eqnarray}
&&\left(\frac{5}{3}\right)^{n_m}\left(\frac{1+m_0(y_{m-1})}{m_0(y_{m-1})}g_x(z_0,z_0)-g_x(z'_0,z_0)\right)\phi_{z_0}^{n_m}(t)\label{eq:5.16}\\
&&+\left(\frac{5}{3}\right)^{n_m}\left(\frac{1+m_0(y_{m-1})}{m_0(y_{m-1})}g_x(z_0,z'_0)-g_x(z'_0,z'_0)\right)\phi_{z'_0}^{n_m}(t).\nonumber
\end{eqnarray}
By substituting the value of $g_x(z_0,z_0)=g_x(z'_0,z'_0)=\frac{m_0(y_{m-1})+m_0(y_{m-1})^2}{2m_0(y_{m-1})+1}\left(\frac{3}{5}\right)^{n_m}$  and
$g_x(z_0,z'_0)=g_x(z'_0,z_0)=\frac{m_0(y_{m-1})^2}{2m_0(y_{m-1})+1}\left(\frac{3}{5}\right)^{n_m}$ into $\eqref{eq:5.16}$, it is easy to
verify that
\begin{equation}
\left(\frac{5}{3}\right)^{n_m}\left(\frac{1+m_0(y_{m-1})}{m_0(y_{m-1})}g_x(z_0,z_0)-g_x(z'_0,z_0)\right)=1,\label{eq:5.17}
\end{equation}
and
\begin{equation}
\left(\frac{5}{3}\right)^{n_m}\left(\frac{1+m_0(y_{m-1})}{m_0(y_{m-1})}g_x(z_0,z'_0)-g_x(z'_0,z'_0)\right)=0.\label{eq:5.18}
\end{equation}
So the coefficient of $v(z_0)$ is $\phi_{z_0}^{n_m}(t)$.

Next we consider the general case. Suppose
$z_0=\tilde{F}_{\omega}(q_0)\in T_x^l$ with $1\leq l<m$. We need to
compute the coefficient of $v(z_0)$. The previous discussion
immediately shows that the contribution of terms in $G_{\Omega_x}^l$
to $v(z_0)$ is $\phi_{z_0}^{n_l}(t)$. Now we consider the terms in
$G_{\Omega_x}^m-G_{\Omega_x}^l$. Let $z_1=\tilde{F}_{\omega 1}(q_0)$
and $z_2=\tilde{F}_{\omega 2}(q_0)$. Notice that in all the terms in
$G_{\Omega_x}^m-G_{\Omega_x}^l$ that contribute to $v(z_0)$, only
those which contain $\phi_{z_1}^{n_{l+1}}(s)$ or
$\phi_{z_2}^{n_{l+1}}(s)$ have supports intersecting the boundary
$S(x)$. Moreover, in calculating the energy
$\mathcal{E}_{\Omega_x}(\phi_{z_i}^{n_{l+1}},v)$, only the part
$\phi_{z_i}^{n_{l+1}}|_{\tilde{\tilde{F}}_{\omega i}(\mathcal{SG})}$
is involved in contributing to the coefficient of $v(z_0)$,
 for i=1,2. Comparing to the standard $\mathcal{SG}$ case, the
 function $\phi_{z_i}^{n_{l+1}}(s)$ has been redefined, but the
 restriction of it to $\tilde{\tilde{F}}_{\omega i}(\mathcal{SG})$ is unchanged.
 So the total contribution of $G_{\Omega_x}^m-G_{\Omega_x}^l$ to
 $v(z_0)$ is as same as the standard case, namely
 $\phi_{z_0}^{n_m}(t)-\phi_{z_0}^{n_l}(t)$. Thus we get that in
 $\mathcal{E}_{\Omega_x}(G_{\Omega_x}^m,v)$, the coefficient of
 $v(z_0)$ is
 $\phi_{z_0}^{n_m}(t)$,
 as required.

 Thus we have proved $\eqref{eq:5.11}$. $\Box$

 \textbf{Theorem 5.4.} \textit{For continuous $F$, the normal derivative of the solution $u$ given by $\eqref{eq:5.10}$ is continuous on $S(x)$.}

 \textit{Proof.} From Theorem 5.3,
 \begin{equation}
 \partial_n u|_{S(x)}=\sum_{m\geq 1}\sum_{z,z'\in
 T_x^m}g_x(z,z')\partial_n\phi_z^{n_m}|_{S(x)}\int_{\Omega_x}\phi_{z'}^{n_m}(t)F(t)dt,\label{eq:5.19}
 \end{equation}
since only those terms containing $\phi_z^k$ whose supports
intersect $S(x)$ contribute to the value of $\partial_n
u|_{S(x)}$.

For fixed $m$, let $z=\tilde{F}_\omega(q_0)\in T_x^m$. Note that
on the cell $\tilde{F}_\omega(\mathcal{SG})$,
$\phi_z^{n_m}=h_0^{y_m}\circ \tilde{F}_{\omega}^{-1}$. By Lemma
3.1, we have
\begin{equation}
\partial_n\phi_z^{n_m}|_{S(x)}=-2\left(\frac{5}{3}\right)^{n_{m+1}}(1-m_0(y_m))2^m\chi_{S_{\omega}(x)}.\label{eq:5.20}
\end{equation}
On the other hand, for $z,z'\in T_x^m$, $g_x(z,z')$ is bounded above
by multiplies of $m_0(y_{m-1})\left(\frac{3}{5}\right)^{n_m}$, hence
by multiplies of $\left(\frac{3}{5}\right)^{n_{m+1}}$ using
$\eqref{eq:2.7}$. It is also easy to see that
$\int_{\Omega_x}\phi_{z'}^{n_m}(t)F(t)dt$ is bounded above by
multiplies of $\frac{1}{3^{n_m}}\parallel F\parallel_{\infty}$.
Combing these estimates with $\eqref{eq:5.20}$, we conclude that
$|\partial_n u||_{S(x)}$ is bounded above by multiplies of
\begin{equation}
\sum_{m\geq1}\sum_{|\omega|=m}\frac{2^m}{3^{n_m}}\parallel
F\parallel_\infty\chi_{S_\omega(x)}.\label{eq:5.21}
\end{equation}
From $\eqref{eq:5.21}$, one can easily verify that $\partial_n u$ is
continuous on $S(x). \Box$

\section*{Acknowledgement} We are grateful to Rachel Kogan who contributed to the proofs of the results in section 2 and 3.

\end{document}